\documentclass{article}

\usepackage[lined]{algorithm2e}
\usepackage{amsmath, amssymb}
\usepackage{makeidx, graphicx, enumerate}

\begin{document}

\title{Eigenvalues and Eigenvectors of the Matrix of Permutation Counts}
\author{Pawan Aurora and Shashank K Mehta\\
Indian Institute of Technology, Kanpur - 208016, India\\
{\small paurora@cse.iitk.ac.in, skmehta@cse.iitk.ac.in}}

\date{}
\maketitle

\newtheorem{definition}{Definition}
\newtheorem{observation}{Observation}
\newtheorem{claim}{Claim}
\newtheorem{problem}{Problem}
\newtheorem{corollary}{Corollary}
\newtheorem{lemma}{Lemma}
\newtheorem{theorem}{Theorem}

\newcommand{\proof}[1]{\noindent \textbf{\textit{Proof}} #1
$\qquad \blacksquare$}
\newcommand{\proofs}[1]{\noindent \textbf{\textit{Proof
Sketch}} #1 $\qquad \blacksquare$}
\newcommand{\ar}{\rightarrow}

\section{Matrix $B$}
Define a $(n^4+n^2)/2\times (n^4+n^2)/2$ symmetric $B$. $(ij)(kl)$ is an index where $i,j,k,l\in [n]$,
$(ab)$ is an unordered pair and $(kl)$ is an ordered pair when $i\neq j$, otherwise it is also an unordered
pair. $B((ij)(kl),(ab)(xy))$ is equal to the number of permutations of $S_n$ in which $\min\{i,j\}$ maps to $k$, 
$\max\{i,j\}$ maps to $l$, $\min\{a,b\}$ maps to $x$ and $\max\{a,b\}$ maps 
to $y$. Clearly each entry of $B$ can be $0$,$(n-4)!$, $(n-3)!$, $(n-2)!$, or $(n-1)!$. 
It can be shown that $B$ is positive semidefinite.

In the following section we will show that $B$ has four distinct eigenvalues: $(3/2)n!$, $n(n-3)!$, 
$(n-1)!/(n-3)$, $2n(n-2)!$ and the corresponding eigenspace dimensions are $1$, ${{n-1}\choose{2}}^2$, 
$({{n-1}\choose{2}}-1)^2$, $(n-1)^2$ respectively.

\section{First Eigenvalue}

The first eigenvalue has only one eigenvector $v$ which is given by $v_{(aa)(xx)} = n-1$ for 
all $a,x$ and $v_{(ab)(xy)}=1$ for all $a\neq b$ and all $x\neq y$. All other entries of $v$
are zero. Let $B_{(ij)(kl)}$ denote the row of $B$ with index $(ij)(kl)$. From the definition
of $B$ it is obvious that row vectors
$B_{(ij)(kl)}$ are non-zero only if either $i=j$ and $k=l$ or $i\neq j$ and $k\neq l$. We will 
compute the inner-product of $B_{(ij)(kl)}$ and $v$ for these two cases.

In inner product $B_{(ii)(kk)}\cdot v$ the details of various terms are as follows. 

\begin{tabular}{llll}
Term Index   &    Term Value   &       No. of terms & total contribution\\
$(ii)(kk)$   &  $(n-1)(n-1)!$  &       $1$                  &  $(n-1)(n-1)!$\\
$(ib)(ky)$   &  $1(n-2)!$      &       $(n-1)^2$            &  $(n-1)^2(n-2)!$ \\
$(aa)(xx)$   &  $(n-1)(n-2)!$  &       $(n-1)^2$            &  $(n-1)^2(n-2)!$\\
$(ab)(xy)$   &  $1(n-3)!$      &       $2{{n-1}\choose{2}}^2$ &  $2(n-3)!{{n-1}\choose{2}}^2$
\end{tabular}

So the inner-product is $(n-1)(n-1)! + (n-1)^2(n-2)! + (n-1)^2(n-1)! + {{n}\choose{2}}(n-1)!
= (n-1)3n!/2 = v_{(ii)(kk)}3n!/2$.

Same details in $B_{(ij)(kl)}\cdot v$ are as follows. 

\begin{tabular}{llll}
Term Index   &     Value       &    No. of terms   & total contribution\\
$(ij)(kl)$   &  $1(n-2)!$      &    $1$                    & $(n-2)!$\\
$(ii)(kk)$   &  $(n-1)(n-2)!$  &    $1$                    & $(n-1)(n-2)!$\\
$(jj)(ll)$   &  $(n-1)(n-2)!$  &    $1$                    & $(n-1)(n-2)!$\\
$(ib)(ky)$   &  $(n-3)!$       &    $(n-2)^2$              & $(n-2)^2(n-3)!$\\
$(jb)(ly)$   &  $(n-3)!$       &    $(n-2)^2$              & $(n-2)^2(n-3)!$\\
$(aa)(xx)$   &  $(n-1)(n-3)!$  &    $(n-2)^2$              & $(n-1)(n-2)^2(n-3)!$\\
$(ab)(xy)$   &  $1(n-4)!$      &    $2{{n-2}\choose{2}}^2$ & $2(n-4)!{{n-2}\choose{2}}^2$
\end{tabular}

In this case the inner-product is $(n-2)!+2(n-1)! + 2(n-2)(n-2)! + (n-2)(n-1)! + 
(n-2)!{{n-2}\choose{2}}$. This simplifies to $3n!/2 = v_{(ij)(kl)}3n!/2$. Hence
$v$ is an eigenvector with eigenvalue $3n!/2$.

\section{Second Eigenvalue}
Let $G_1=(V_1,E_1,w_1)$ be an undirected and $G_2=(V_2,E_2,w_2)$ be a directed 
edge-weighted graphs, where $V_1$ and $V_2$ are subsets of $[n]$. 
Neither graph has loop-edges and $(x,y)\in E_2$ if and only if $(y,x)\in E_2$.
For each $e_1=(ab)\in E_1$ and $e_2=(xy)\in E_2$ we associate $(e_1,e_2)$ with axis
$(ab,xy)$. Note $(ab,xy)$ is different from $(ab,yx)$. We define a vector $v(G_1,G_2)$ 
or simply $v$ as follows: if $(ab)\in E_1$ and $(xy)\in E_2$ then 
$v_{(ab)(xy)} = w_1(ab).w_2(xy)$, otherwise $v_{(ab)(xy)} = 0$. In this subsection
and in the following subsection we will show $v(G_1,G_2)$ are eigenvectors for
various instances of $G_1$ and $G_2$.

\begin{figure}
\begin{center}
\includegraphics[scale=0.5]{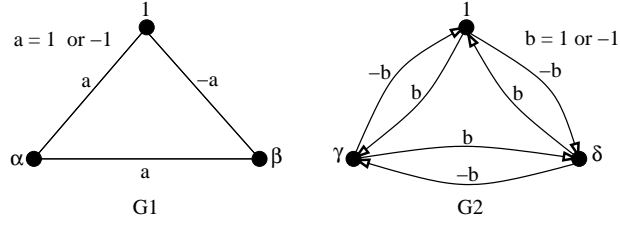}
\end{center}
\label{eigen2}
\caption{Index Graphs for the Eigenvectors for Second Eigenvalue}
\end{figure}

We will use following notations regarding $G_1$ and $G_2$. $\delta(a)$ will denote the
edges incident on vertex $a\in V_1$ and $N(a)$ will denote the neighbors of vertex $a\in V_1$.
$N[a]$ will denote the closed neighborhood of $a$ which is $N(a)\cup \{a\}$.
Graph $G_2$ is directed so $\delta^{\rightarrow}(x)$ denotes the set of outgoing edges
and $\delta^{\leftarrow}(x)$ denotes the incoming edges incident on $x\in V_2$. Similarly
neighborhood will be denoted by $N^{\rightarrow}(x)$ and $N^{\leftarrow}(x)$ respectively.

In this subsection we assume following properties of the weight-functions:
(i) all weights are $1$ or $-1$; (ii) $\sum_{b\in N(a),a<b}w_1(a,b) -\sum_{b\in N(a),a>b}w_1(a,b) 
= 0$; (iii) $w_2$ be such that $w_2(x,y)=-w_2(y,x)$ for all $(x,y)\in E_2$,
and (iv) $\sum_{y\in N^{\rightarrow}[x]}w_2(x,y) = \sum_{y\in N^{\leftarrow}[x]}w_2(y,x) = 0$
for each $x\in V_2$.

\begin{lemma} \label{lem3}The inner-product $B_{(ij)(kl)}\cdot v(G_1,G_2)$
is zero unless $i,j\in V_1$ and $k,l\in V_2$.
\end{lemma}

\proof{
If $k,l$ both are out of $V_2$ then for each $(a,b)$ we will get contributions from $(x,y)$ and
$(y,x)$ which will cancel each other.

Next suppose $k\in V_2$ and $l\notin V_2$. Here we consider two cases: $i,j\notin V_1$ and 
$i\in V_1,j\notin V_1$. In the first case for every non-zero term with index $(a,b)(x,y)$, 
both $x$ and $y$ will be different from $k$ because $k$ is the image of $i$ so neither $x$ 
nor $y$ can be image of $i$. Hence
contributions from $(x,y)$ and $(y,x)$ will cancel each other. If $i\in V_1$
and $a$ and $b$ are different from $i$, then the situation will be same as above. 
On the other hand if $a=i$, then we will get a total sum as $\sum_{y\in N^{\rightarrow}(k)}
w_1(a,b)w_2(x,y)$ if $a<b$. In case $a>b$, then the total will be 
$\sum_{y\in N^{\rightarrow}(k)} w_1(a,b)w_2(y,x)$. In both cases the sum is zero.

Next consider the case when $i\in V_1, j\notin V_1,k\in V_2, l\in V_2$. 
If $a,b$ are both different from $i$, then
contribution of $(x,y)$ and of $(y,x)$ will cancel each other.
Now, let $a=i$. In this case also we consider two sub-cases: $(kl)\notin E_2$
and $(k,l)\in E_2$. In the first sub case the contribution will be 
$(n-3)!\sum_{y\in N^{\rightarrow}(k)}w_1(i,b)w_2(k,y)$ or
$(n-3)!\sum_{y\in N^{\rightarrow}(k)} w_1(i,b)w_2(y,k)$ depending on $i<b$ or $i>b$. 
In both cases
the sum is zero. In the second sub case case, assume $i<b$. Then the sum of the contributions
is $\sum_{y\in N^{\rightarrow}(k)\setminus \{l\}} w_1(i,b)w_2(k,y)$. This simplifies to
$-w_1(i,b)w_2(k,l)$. If $i>b$, then we get $-w_1(i,b)w_2(l,k)=
w_1(i,b)w_2(k,l)$. In this case if we add up the contributions for all $b\in N(i)$,
then we get $(n-3)!(\sum_{b\in N(i): b<i} w_1(i,b)w_2(k,l) - 
\sum_{b\in N(i): b>i} w_1(i,b)w_2(k,l))$. This reduces to zero because it is given that
$\sum_{b\in N(i): b<i} w_1(i,b) - \sum_{b\in N(i): b>i} w_1(i,b)=0$.

Finally let $i\in V_1, j\in V_1,k\in V_2, l\notin V_2$. If $a\neq i$, then
contribution of $(x,y)$ will cancel that of $(y,x)$. If $a=i$, then the
contribution will be
$(n-3)!\sum_{y\in N^{\rightarrow}(k)}w_1(i,b)w_2(k,y)$ or
$(n-3)!\sum_{y\in N^{\rightarrow}(k)} w_1(i,b)w_2(y,k)$ depending on $i<b$ or $i>b$. 
As observed earlier, in both cases the sum is zero.
}

\begin{theorem} For graphs $G_1$ and $G_2$ given in figure \ref{eigen2} where $\alpha < \beta$
and $\gamma < \delta$, $v(G_1,G_2)$ is an eigenvector with eigen value $n(n-3)!$.
\end{theorem}

\proof{
We have seen from Lemma \ref{lem3} that in the inner-product $B_{(ij)(kl)}\cdot v$ the only 
terms, $(ab)(xy)$, have non-zero value are those in which $a,b \in V_1$ and $x,y\in V_2$.
Without loss of generality assume that $i<j$. Denote the third vertex in $V_1$ apart
from $i,j$ by $p$ and the third vertex in $V_2$ apart from $k,l$ by $q$.

The non-zero contribution is possible for three terms: (i) $(ab)(xy) = (ij)(kl)$,
(ii) $(ab)(xy) = (ip)(kq)$, (iii) $(ab)(xy) = (jp)(lq)$. The values of these terms are as follows:

(i) $(n-2)!w_1(ij)w_2(kl)$;

(ii) if $i<p$, then $(n-3)!w_1(ip)w_2(kq)$, otherwise $(n-3)!w_1(ip)w_2(qk)$;

(iii) if $j<p$, then $(n-3)!w_1(jp)w_2(lq)$, otherwise $(n-3)!w_1(jp)w_2(ql)$.

The sum of these terms is the inner-product. We determine it for all the three cases:

Case $p<i<j$: Recall that for $V_1=\{L,M,H\}$ with $L<M<H$, properties of $w_1$ implies $w_1(LM)-w_1(MH)=0$,
$w_1(LM) + w_1(LH)=0$, and $w_1(LH)+w_1(MH)=0$. So the total is $(n-2)!w_1(ij)w_2(kl) + 
(n-3)!(w_1(ip)w_2(qk)-w_1(jp)w_2(ql))$. The $w_1$-equations given above simplifies the expression
to $(n-2)!w_1(ij)w_2(kl) + (n-3)!(w_1(ij)w_2(kl)-w_1(jp)(-w_2(kl)))$. This expression simplifies to
$n(n-3)!w_1(ij)w_2(kl)$.

Case $i<p<j$: Here the total is $(n-2)!w_1(ij)w_2(kl) + (n-3)!(w_1(ip)w_2(kq)+ w_1(jp)w_2(ql))$
which is equal to $(n-2)!w_1(ij)w_2(kl) + (n-3)!((-w_1(ij))(-w_2(kl))\\
+(-w_1(jp))(-w_2(kl)))$.
This also simplifies to $n(n-3)!w_1(ij)w_2(kl)$.

Case $i<j<p$: In this case total is $(n-2)!w_1(ij)w_2(kl) + (n-3)!(w_1(ip)w_2(kq)+ w_1(jp)w_2(lq))$
which is equal to $(n-2)!w_1(ij)w_2(kl) + (n-3)!((-w_1(ij))(-w_2(kl))+w_1(jp)w_2(kl))$.
This too simplifies to $n(n-3)!w_1(ij)w_2(kl)$.

Hence in each the inner-product is $n(n-3)!v[(ij)(kl)]$. So $B\cdot v = n(n-3)!v$, i.e.,
$v$ is an eigenvector with eigenvalue $n(n-3)!$.
}

\begin{theorem} The eigenspace of $B$ corresponding to eigenvalue $n(n-3)!$ has
dimension at least ${{n-1}\choose{2}}^2$.
\end{theorem}

\proof{We can define one vector $v(G_1,G_2)$ for each choice of $(\alpha,\beta)$
and each choice of $(\gamma,\delta)$ from $2,\dots,n$. Let $v$ and $v'$ be 
such a vector due to $(\alpha,\beta,\gamma,\delta)$.
Then $v_{(\alpha,\beta),(\gamma,\delta)}$ is non-zero but the same component 
of all other eigenvectors of this type is zero. Hence these vectors are independent.
}

\section{Third Eigenvalue}

The eigenvectors are defined in the same way as in the previous section.
In this subsection we assume following properties of the weight-functions associated
with $G_1$ and $G_2$: (i) all weights are $1$ or $-1$; (ii) $\sum_{b\in N(a)}w_1(a,b)=0$ 
for all $a\in V_1$; (iii) $w_2(x,y)=w_2(y,x)$ for all $(xy)\in E_2$; and (iv) 
for all $x$, $\sum_{y\in N^{\rightarrow}}(x) w_2(x,y)= \sum_{y\in N^{\leftarrow}}(x) w_2(y,x) = 0$.

\begin{lemma}\label{lem4} The inner product $B_{(ij)(kl)}\cdot v(G_1,G_2)$ is non-zero only if
$(ij)\in E_1$ and $(kl)\in E_2$.
\end{lemma}

\proof{

{\bf Case} $k\notin V_2$, $l\notin V_2$\\
Clearly if $a\in \{i,j\}$ or $b\in \{i,j\}$, then $B_{(ij)(kl)}[(ab)(xy)]=0$. Otherwise
$\sum_{(xy)\in E_2}B_{(ij)(kl)}[(ab)(xy)].v[(ab)(xy)] = (n-4)!\sum_{x\in V_2}
\sum_{y\in N^{\rightarrow}(x)}v[(ab)(xy)] = 0$. So $B_{(ij)(kl)}\cdot v = 0$.

{\bf Case} $i\notin V_1$, $j\notin V_1$\\
As in the above case $B_{(ij)(kl)}\cdot v = 0$.

{\bf Case} $i\in V_1$, $j\notin V_1$\\
Assume that $i<j$.

{\em Sub-case} $x\neq k$\\
First let $a\notin N[i]$.

$\sum_{b\in N[a]}B_{(ij)(kl)}[(ab)(xy)].v[(ab)(xy)]
= (n-4)!\sum_{b\in N[a]}v[(ab)(xy)] = 0$.

Next $a\in N(i)$.

$\sum_{a\in N(i)}\sum_{b\in N(a)}B_{(ij)(kl)}[(ab)(xy)].v[(ab)(xy)]\\
= (n-4)!\sum_{a\in N(i)}(-v[(ai)(xy)]) = 0$.

{\em Sub-case} $x=k$\\

$\sum_{b\in N(i)}B_{(ij)(kl)}[(ab)(xy)].v[(ab)(ky)]
= (n-3)!\sum_{b\in N(i)}v[(ab)(ky)] = 0$.

Similar argument works for $j<i$.

{\bf Case} $k\in V_2$, $l\notin V_2$ or $k\notin V_2$, $l\in V_2$\\
This case is similar to the previous case.

At this stage in the proof we have shown that the inner product can be non-zero only 
if $i\in V_1,j\in V_1, k\in V_2$ and $l\in V_2$. Now we will show that even in these 
cases the inner product will be zero if $(ij)\notin E_1$ or $(kl)\in E_2$.

For the remainder of the proof $i\in V_1,j\in V_1, k\in V_2$ and $l\in V_2$.
For the remaining three cases we assume that $(i,j)\notin E_1$.

{\em Case} $a \notin N[i]\cup N[j]$\\
Let $i<j$.

$\sum_{b\in N(a)}B_{(ij)(kl)}[(ab)(xy)].v[(ab)(xy)]
= (n-4)!\sum_{b\in N(a)}v[(ab)(xy)] = 0$.

{\bf Case} $a\in (N(i)\cup N(j))\setminus \{i,j\}$\\

Recall that $i$ and $j$ are not adjacent.

$\sum_{b\in N(a)\setminus \{i\}}B_{(ij)(kl)}[(ab)(xy)].v[(ab)(xy)]$ is equal to 
$(n-4)!(-v[(ai)(xy)])$ if $a\in N(i)\setminus N(j)$; $(n-4)!(-v[(aj)(xy)])$ if $a\in N(j)\setminus N(i)$;
and $(n-4)!(-v[(ai)(xy)]-v[(aj)(xy)])$ if $a\in N(i)\cap N(j)$;

So $\sum_{a\in (N(i)\cup N(j))\setminus \{i,j\}}\sum_{b\in N(a)\setminus \{i\}}B_{(ij)(kl)}[(ab)(xy)].v[(ab)(xy)]$
is equal to $|N(i)\cap N(j)|(n-4)!\sum_{a\in N(i)}(-v[(ai)(xy)]) + (n-4)!\sum_{a\in N(j)}(-v[(aj)(xy)]) = 0$.

{\bf Case} $a=i$ or $a=j$\\
$\sum_{b\in N(i)}B_{(ij)(kl)}[(ab)(xy)]v[(ab)(xy)] = (n-3)! \sum_{b\in N(i)}v[(ib)(xy)] = 0$.
Similarly for $a=j$. 

From the last three cases we conclude that if $(ij)$ is not an edge in $G_1$ then 
$B_{(ij)(kl)}\cdot v(G_1,G_2)=0$. Similarly we can show that $B_{(ij)(kl)}\cdot v(G_1,G_2)=0$
if $(kl)$ is not an edge in $G_2$.
}

\begin{figure}
\begin{center}
\includegraphics[scale=0.5]{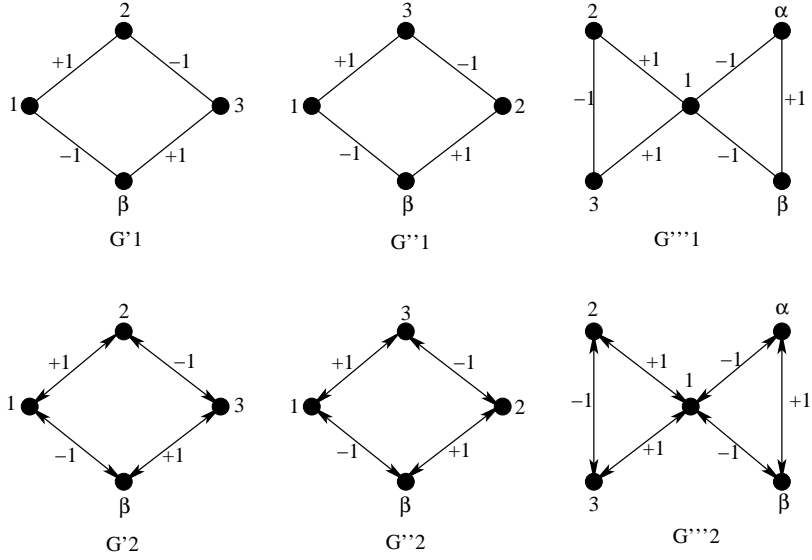}
\end{center}
\label{eigen3}
\caption{Index Graphs for the Eigenvectors for Third Eigenvalue}
\end{figure}

In this case $G_1$ is any graph among $G'_1,G''_1,$ and $G'''_1$ and $G_2$ is any graph
among $G'_2,G''_2,$ and $G'''_2$. These graphs are given in figure \ref{eigen3}.

\begin{theorem} $v(G_1,G_2)$ is an eigenvector of $B$ for each $G_1\in \{G'_1,G''_1,G'''_1\}$
and $G_2\in \{G'_2,G''_2,G'''_2\}$ with eigenvalue $(n-1)!/(n-3)$.
\end{theorem}

\proof{
From lemma \ref{lem4} we have seen that $B_{(ij)(kl)}\cdot v(G_1,G_2)$ is zero if $(ij)\notin E_1$ or
$(kl)\notin E_2$. For $(ij)\in E_1$ and $(kl)\in E_2$ there are $2\times 14\times 14$ cases. 
Each of these cases is equivalent to one of the following 15 cases: 
(1) $B_{(12)(12)}\cdot v(G'_1,G'_2)$, 
(2) $B_{(12)(23)}\cdot v(G'_1,G'_2)$, 
(3) $B_{(12)(12)}\cdot v(G'_1,G'''_2)$, 
(4) $B_{(12)(21)}\cdot v(G'_1,G'''_2)$, 
(5) $B_{(12)(23)}\cdot v(G'_1,G'''_2)$, 
(6) $B_{(12)(12)}\cdot v(G'''_1,G'_2)$, 
(7) $B_{(12)(23)}\cdot v(G'''_1,G'_2)$, 
(8) $B_{(12)(12)}\cdot v(G'''_1,G'''_2)$, 
(9) $B_{(12)(21)}\cdot v(G'''_1,G'''_2)$, 
(10) $B_{(12)(23)}\cdot v(G'''_1,G'''_2)$, 
(11) $B_{(23)(12)}\cdot v(G'''_1,G'_2)$, 
(12) $B_{(23)(23)}\cdot v(G'''_1,G'_2)$, 
(13) $B_{(23)(12)}\cdot v(G'''_1,G'''_2)$, 
(14) $B_{(23)(21)}\cdot v(G'''_1,G'''_2)$. 
(15) $B_{(23)(23)}\cdot v(G'''_1,G'''_2)$. 

The theorem can be established by showing $B_{(ij)(kl)}\cdot v = ((n-1)!/(n-3))v[(ij)(kl)]$
for each of the cases enumerated above. We illustrate the same for cases 1 and 13.

Case (1): $B_{(12)(12)}\cdot v[G'_1,G'_2]$ is equal to $B_{(12)(12)}[(12)(12)]v[(12)(12)]
+ B_{(12)(12)}[(1\beta)(1\delta)]v[(1\beta)(1\delta)] + B_{(12)(12)}[(23)(23)]v[(23)(23)] +\\ 
B_{(12)(12)}[(3\beta)(3\delta)]v[(3\beta)(3\delta)] + B_{(12)(12)}[(3\beta)(\delta 3)]$.
This is equal to $(n-2)! + 2(n-3)! + 2(n-4)! = (n-1)!/(n-3)$.
Since $v[(12)(12)=1$, it is also equal to $((n-1)!/(n-3))v[(12)(12)]$.

Case (10): $B_{(12)(23)}\cdot v(G'''_1,G'''_2)$ is equal to 
$B_{(12),(23)}[(12),(23)]v[(12),(23)]
+B_{(12),(23)}[(13),(21)]v[(13),(21)]
+B_{(12),(23)}[(1\alpha),(21)]v[(1\alpha),(21)]\\
+B_{(12),(23)}[(1\beta),(21)]v[(1\beta),(21)]
+B_{(12),(23)}[(23),(31)]v[(23),(31)]\\
+B_{(12),(23)}[(\alpha\beta),(1\gamma)]v[(\alpha\beta),(1\gamma)]
+B_{(12),(23)}[(\alpha\beta),(\gamma 1)]v[(\alpha\beta),(\gamma 1)]\\
+B_{(12),(23)}[(\alpha\beta),(1\delta)]v[(\alpha\beta),(1\delta)]
+B_{(12),(23)}[(\alpha\beta),(\delta 1)]v[(\alpha\beta),(\delta 1)]\\
+B_{(12),(23)}[(\alpha\beta),(\gamma\delta)]v[(\alpha\beta),(\gamma\delta)]
+B_{(12),(23)}[(\alpha\beta),(\delta\gamma)]v[(\alpha\beta),(\delta\gamma)]$
This expression is equal to $-(n-2)!-2(n-3)!-2(n-4)!$ which is equal to $-(n-1)!/(n-3)$.
Since $v[(12)(23)]=-1$, $B_{(12)(23)}\cdot v(G'''_1,G'''_2)=((n-1)!/(n-3))v[(12)(23)]$.
}

\begin{theorem} The dimension of the eigenspace of eigenvalue $(n-1)!/(n-3)$ is $({n-1}\choose{2} - 1)^2$.
\end{theorem}

\proof{Each vector described above has an entry which is non-zero while the same entry is zero among
the other vectors. In $v[G'_1(\beta),G'_2(\delta)]$ such an entry is $(3\beta)(3\delta)$; in 
$v[G'_1(\beta),G''_2(\delta)]$ such an entry is $(3\beta)(2\delta)$;in
$v[G'_1(\beta),G'''_2(\gamma\delta)]$ such an entry is $(3\beta)(\gamma\delta)$ so on. Hence all these 
vectors are linearly independent. Therefore the dimension of this eigenspace is at least equal to
the number of these vectors. The number of instances $G'_1$ are $n-3$, that of $G''_1$ is $n-3$ and that 
of $G'''_1$ is ${n-3}\choose{2}$. So the total number of instances of $G_1$ is ${n-1}\choose{2}-1$.
The number of instances of $G_2$ is the same. Hence the total number of vectors is $({n-1}\choose{2}-1)^2$.}

\section{Fourth Eigenvalue}

To define the eigenvectors for this eigenvalue consider three $n\times n$ matrices:
$A_1,A_2,A_3$, given below. The symmetric vectorization $v(E)$ of block-matrix $E$
will be shown to an eigenvector of $B$ with eigenvalue $2n(n-2)!$.

\[ A_1=\left[ \begin{array}{cccccc}
1-(n-1)(n-2) & 0 & \dots & 0 & 0 & 0\\
0            & n-1 & \dots & 0 & 0 & 0\\
\dots        & \dots & \dots & \dots & \dots & \dots\\
0            & 0     & \dots & 0 & n-1 & 0\\
0            & 0     & \dots & 0 & 0 & -1 \end{array} 
\right ] \]

\[ A_2=\left [ \begin{array}{cccccc}
0       & -n+2+\frac{1}{n-2} & -n+2+\frac{1}{n-2} & \dots       & -n+2+\frac{1}{n-2} & -n+2\\
\frac{1}{n-2} & 0              & \frac{n-1}{n-2}      & \dots       & \frac{n-1}{n-2}    & 1\\
\frac{1}{n-2} & \frac{n-1}{n-2}    & 0                & \dots       & \frac{n-1}{n-2}    & 1\\
\dots   & \dots          & \dots                      & \dots       & \dots           & \dots\\
\frac{1}{n-2} & \frac{n-1}{n-2}    & \frac{n-1}{n-2}  & \dots       &  0             & 1\\
-1      & -1             & -1                         & \dots       &  -1            & 0 \end{array}
\right ] \]


\[ A_3=\left [ \begin{array}{ccccccc}
0       & -n+2 & -n+2 & -n+2 & \dots       & -n+2 & -n+3\\
n-2   & 0      & 0      & 0      & \dots       & 0       &  1\\
n-2   & 0      & 0      & 0      & \dots       & 0       &  1\\
\dots   & \dots  & \dots  & \dots  & \dots       & \dots   & \dots\\
n-2   & 0      & 0      & 0      & \dots       & 0       &  1\\
n-3   & -1     & -1     & -1     & \dots       & -1      & 0 \end{array}
\right ] \]

\[E= \left[ \begin{array}{ccccc}
A_1 & A_3 & A_2 & \dots & A_2\\
A_3^T & -A_1 & -A_2 & \dots & -A_2\\
A_2^T & -A_2^T & 0 & \dots & 0\\
\dots & \dots & \dots & \dots & \dots\\
A_2^T & -A_2^T & 0 & \dots & 0  \end{array} \right ]
\]

Let $X(a_1\ar a'_1,a_2\ar a'_2,a_3 \ar a'_3,a_4\ar a'_4)$ denote the number of permutations in which 
$a_i$ maps to $a'_i$ for $i=1,2,3,4$. Since matrix $A_1$ is diagonal, 
$(B\cdot v)_{ik,jl}$ can be expressed as operations on the 
$n\times n$ blocks as follows: $(B\cdot v)_{ik,jl} = \sum_{pr,qs} B_{ik,jl}(pr,qs)v(pr,qs)
= \sum_{q,s}A_1(q,s).X(i\ar j,k\ar l,1\ar q,1\ar s)
- \sum_{q,s}A_1(q,s).X(i\ar j,k\ar l,2\ar q,2\ar s)
+ \sum_{q,s}A_3(q,s).X(i\ar j,k\ar l,1\ar q,2\ar s)
+ \sum_{r>2}\sum_{q,s}A_2(q,s).X(i\ar j,k\ar l,1\ar q,r\ar s)
- \sum_{r>2}\sum_{q,s}A_2(q,s).X(i\ar j,k\ar l,2\ar q,r\ar s) $.
Denote this expression by $C(ik)_{jl}$.
Now we will evaluate the value of this expression for various values of $(ik,jl)$.

Case: $ik=11$. In this case we need to show that $C(11) = 2n(n-2)!A_1$.

For $jl=11$ we have $C(11)_{11} = 
(n-1)!(1-(n-1)(n-2)) - (-1+(n-1)(n-2))(n-2)! + (-(n-2)^2-(n-3))(n-2)! 
(n-2)((n-2)((-n+2+(1/(n-2)))) - (n-2))(n-2)! - (n-2)^2(1+ (n-3)(n-1)/(n-2))(n-3)!$. This expression
simplifies to $2n(n-2)!(1-(n-1)(n-2)) = 2n(n-2)!A_1(1,1)$.

For $jl=22$ we have $C(11)_{22} = 
(n-1)(n-1)! 
- ((1-(n-1)(n-2)) + (n-1)(n-3)-1)(n-2)!
+((n-2)+1)(n-2)! 
+ (n-2)(1/(n-2) +(n-3)(n-1)/(n-2) +1)(n-2)! 
- (n-2)((n-3)(-n+2-1/(n-2)) - (n-2) + (1/(n-2) + (n-4)(n-1)/(n-2) +1)(n-3) - 1)(n-3)!$.
It simplifies to $2n! = (n-1).2n(n-2)! = 2n(n-2)!.A_1(2,2)$.

For $jl=tt$ where $2\leq t\leq n-1$ we get the similar expression, i.e., $2n(n-2)!.A_1(t,t)$.

In case $jl=nn$, $C(11)_{nn} =
-(n-1)! - (1-(n-1)(n-2) + (n-1)(n-2))(n-2)! + ((n-3)-(n-2))(n-2)! -(n-2)(n-2)!
-(n-2)((-n+2+(1/(n-2))(n-2) + ((1/(n-2)) + ((n-3)(n-1)/(n-2))(n-2))(n-3)!
= -2n(n-2)! = 2n(n-2)!A_1(n,n)$.

In every case where $j\neq l$, $X(1\ar j,1\ar l,*,*) = 0$. So $C(11)_{jl}=0$ for all
$j\neq l$. Thus $C(11) = 2n(n-2)!A_1$.


Case: $ik=22$. In this case we have to show that $C(22) = 2n(n-2)!(-A_1)$.

We observe following equalities.
$A_1(q,s).X(2\ar j,2\ar l,2\ar q,2\ar s)=A_1(q,s).X(1\ar j,1\ar l,1\ar q,1\ar s)$,
$A_1(q,s).X(2\ar j,2\ar l,1\ar q,1\ar s)=A_1(q,s).X(1\ar j,1\ar l,2\ar q,2\ar s)$,
$A_2(q,s).X(2\ar j,2\ar l,1\ar q,r\ar s)  =A_2(q,s).X(1\ar j,1\ar l,2\ar q,r\ar s)$,
$A_2(q,s).X(2\ar j,2\ar l,2\ar q,r\ar s)  =A_2(q,s).X(1\ar j,1\ar l,1\ar q,r\ar s)$.
Further,
$A_3(q,s).X(2\ar j,2\ar l,1\ar q,2\ar s) = A_3(q,s).X(1\ar j,1\ar l,2\ar q,1\ar s)
= A_3(s,q).X(1\ar j,1\ar l,1\ar q,2\ar s) = -A_3(q,s).X(1\ar j,1\ar l,1\ar q,2\ar s)$.

Substituting these values in the expression of $C(22)$ we get $C(22) = -C(11)=2n(n-2)!(-A_1)$.

Case $ik=12$. In this case we have to show that $C(12)=2n(n-2)!A_3$.

For $jl=12$, $C(12)_{12} = 
(1-(n-1)(n-2))(n-2)! - (n-1)(n-2)! + ((n-3)(-n+2+(1/(n-2))) - (n-2))(n-2)(n-3)!
- (1+(n-1)(n-3)/(n-2))(n-2)(n-3)! = -2(n-2)n(n-2)! = 2n(n-2)!A_3(1,2)$.

For $jl=21$, $C(12)_{21} =
(n-1)(n-2)! - (1-(n-1)(n-2)) + (n-2)(n-2)! + (n-2)(1+ (n-3)(n-1)/(n-2))(n-3)! -
(n-2)((n-3)(-n+2+1/(n-2))-(n-2))(n-3)! = (n-2)!(2n^2-4n) = 2n(n-2)!A_3(2,1)$.

Next for $jl=1n$, $C(12)_{1n} =
(1-(n-1)(n-2))(n-2)! + (n-2)! - (n-3)(n-2)! - (n-2)^2(n-2-(1/(n-2)))(n-3)! = - 2n(n-3)(n-2)! =
2n(n-2)!A_3(1,n)$.

Expression for $jl=n1$, $C(12)_{n1} = -1.(n-2)!-(1-(n-1)-(n-2))(n-2)! +
(n-3)(n-2)!+0+(n-2)^2(n-2-1/(n-2))(n-3)! = 2n(n-3)(n-2)!=2n(n-2)!A_3(n,1)$.

Expression for $jl=2n$, $C(12)_{2n}=
(n-1)(n-2)! + (n-2)! + (n-2)! +((n-1)(n-3)/(n-2) + 1/(n-2))(n-2)(n-3)! + (n-2)(n-3)!
= 2n(n-2)! = 2n(n-2)!A_3(2,n)$. Similarly it can be seen that $C(12)_{rn} = 2n(n-2)!A_3(r,n)$
for every $2 < r < n-1$.

Expression for $jl=n2$, $C(12)_{n2} = -(n-2)! - (n-1)(n-2)! - (n-2)! + (-1)(n-2)(n-3)!
- (n-2)(1/(n-2) + (n-3)(n-1)/(n-2))(n-3)! = -2n(n-2)! = (-1)A_3(n,2)$. Similarly
it can be shown that $C(12)_{n,r} = (-1)A_3(n,r)$ for $1< r < n-1$.

For $j\in \{2,3,\dots,n-1\}$ and $l\in \{2,3,\dots,n-1\}$, $C(12)_{jl} =
(n-1)(n-2)! - (n-1)(n-2)! + 0 +(1 + 1/(n-2) + (n-1)(n-4)/(n-2))(n-2)(n-3)!
- (1 + 1/(n-2) + (n-1)(n-4)/(n-2))(n-2)(n-3)! = 0$.

Putting these together we see that $C(12) = 2n(n-2)!A_3$.


Case $ik=1k$ with $k\geq 3$: Next we will show that $C(1k) = 2n(n-2)!A_2$, for $k\geq 3$. 
We will show the details of the computation of $C(1,3)$.

For $jl=12$, $C(13)_{12} = (1-(n-1)(n-2))(n-2)! - ((n-3)(n-1)-1)(n-3)! - ((n-3)(n-2)+(n-3))(n-3)!
- (n-2-1/(n-2))(n-2)! + (n-3)(-(n-2-1/(n-2))(n-3)-(n-2))(n-3)! - ((n-3)(n-1)/(n-2))(n-3)!
-(n-3)((n-3)(n-4)(n-1)/(n-2) + (n-3))(n-4)!$.  
It simplifies to $-(n-2)!(2n^2-7n+5)-(n-3)!(3n^2-13n+10)$ which is equal to 
$-(n-2-1/(n-2))2n(n-2)! = 2n(n-2)!A_2(1,2)$. Similarly we can show $C(13)_{1l} = 2n(n-2)!A_2(1,l)$,
for all $2\leq l\leq n-1$.

For $jl=21$, $C(13)_{21} = (n-1)(n-2)! - ((n-3)(n-1)-1)(n-3)! + (n-3)! + (1/(n-2))(n-2)!
+(n-3)(1+(n-3)(n-1)/(n-2))(n-3)! - ((n-3)/(n-2)-1)(n-3)!  - (n-3)(n-3 +(n-3)(n-4)(n-1)/(n-2))(n-4)!
= (1/(n-2)).2n(n-2)! = 2n(n-2)!A_2(2,1)$. Similarly we can show $C(13)_{j1} = 2n(n-2)!A_2(j,1)$,
for all $2\leq j\leq n-1$.

For $jl=1n$, $C(13)_{1n} = (n-2)! - (n-2)(n-1)! - (n-1)! - 2(n-2)(n-2)! - (n-2)(n-3)(n-2)!
+ (n-3)(n-3)! - (n-2)(n-3)! - (n-1)(n-3)(n-3)! = -2(n-2)n(n-2)! = 2n(n-2)!A_2(1,n)$.

For $jl=n1$, $C(13)_{n1} = -(n-2)! - (n-2)(n-1)(n-3)! - (n-2)(n-3)! - (n-2)! + 0 - 1.(n-3)! 
- (n-3)(n-2)((n-3)((n-1)/(n-2))(n-4)! = -2n(n-2)! = 2n(n-2)!A_2(n,1)$.

For $jl=2n$, $C(13)_{2n} = (n-1)(n-2)! - (1-(n-1)(n-2) + (n-3)(n-1))(n-3)! + (n-2)(n-3)! + (n-2)!
+(n-3)(1/(n-2) + (n-3)(n-1)/(n-2))(n-3)! - ((-1)(n-3)! - (n-3)((n-3)/(n-2) - (n-3)(n-2-1/(n-2))
+ (n-3)(n-4)(n-1)/(n-2))(n-3)! = 1.2n(n-2)! = A_2(2n).2n(n-2)!$. Similarly we can show $C(13)_{jn} 
= 2n(n-2)!A_2(j,n)$, for all $2\leq j\leq n-1$.

For $jl=n2$, $C(13)_{n2} = -(n-2)! -(-1-(n-1)(n-2)+(n-1(n-3)))(n-3)! + (n-3-(n-3))(n-3)!
+0-(n-3)(n-3)! - (-(n-2) + 1/(n-2) + (n-3)(n-1)/(n-2))(n-3)! = 0 = 2n(n-2)!.A_2(n,2)$.
Similarly it can be shown that $C(13)_{nl} = 0 = 2n(n-2)!.A_2(n,l)$ for all $2\leq l\leq n-1$.

For $jl=23$, $C(13)_{23} = (n-1)(n-2)! - (1 - (n-1)(n-2) -1 + (n-1)(n-4))(n-3)! - (n-2+1)(n-3)!
+((n-1)/(n-2))(n-2)! + (n-3)((n-3)(n-1)/(n-2))(n-3)! - ((n-4)(n-1)/(n-2) - (n-2) + 1/(n-2))(n-3)!
- ( (n-4)/(n-2) -1 +(n-4) - (n-2) - (n-4)(n-2-1/(n-2)) + (n-4)(n-5)(n-1)/(n-2))(n-3)! = 2n(n-1)(n-2)!/(n-2)
= A_2(2,3).2n(n-2)!$. Similarly we can show that $C(13)_{jl} 
= 2n(n-2)!A_2(j,l)$, for all $2\leq j,l\leq n-1$ and $j\neq l$.

For $jl=n2$, $C(13)_{n2} = (-1)(n-2)! - ((n-1)(n-3)-(n-1)(n-2)+1)(n-3)! + 0 + 0 - (n-3).(n-3)!
-((n-3)(n-1)/(n-2) - (n-2) + 1/(n-2))(n-3)! - (-(n-3)(n-2-1/(n-2)) + (n-3)/(n-2) + (n-3)(n-4)(n-1)/(n-2))(n-3)!=0$.
Similarly we can show $C(13)_{nl} = 2n(n-2)!A_2(n,l)$, for all $2\leq l\leq n-1$.

For $jl=jj$, trivially $C(13)_{jj}=0$. Hence $C(13)_{jj}=2n(n-2)!.A_2(j,j)$.

These observations conclude that $C(13)=2n(n-2)!A_2$. Similarly it can be seen that $C(1,k)=2n(n-2)!A_2$
for all $k\geq 3$.


Case $ik=2k$ for $k\geq 3$: Now we will show that $C(2,k) = 2n(n-2)!(-A_2)$ for all $k\geq 3$. 
We will show the details of the computation of $C(2,3)$.

$C(2,3)_{jl} = \sum_{q,s}A_1(q,s).X(2\ar j,3\ar l,1\ar q,1\ar s)
- \sum_{q,s}A_1(q,s).X(2\ar j,3\ar l,2\ar q,2\ar s)
+ \sum_{q,s}A_3(q,s).X(2\ar j,3\ar l,1\ar q,2\ar s)
+ \sum_{r>2}\sum_{q,s}A_2(q,s).X(2\ar j,3\ar l,1\ar q,r\ar s)
- \sum_{r>2}\sum_{q,s}A_2(q,s).X(2\ar j,3\ar l,2\ar q,r\ar s) $.

We will show that the right hand side is equal to $-C(2,3)_{jl}$.
We state the following facts.

$A_1(q,s).X(2\ar j,3\ar l,1\ar q,1\ar s) = A_1(q,s).X(1\ar j,3\ar l,2\ar q,2\ar s)$.
$A_1(q,s).X(2\ar j,3\ar l,2\ar q,2\ar s)=A_1(q,s).X(1\ar j,3\ar l,1\ar q,1\ar s)a$.
$A_3(q,s).X(2\ar j,3\ar l,1\ar q,2\ar s) A_3(q,s).X(1\ar j,3\ar l,2\ar q,1\ar s)= 
A_3(s,q).X(1\ar j,3\ar l,1\ar q,2\ar s) = -A_3(q,s).X(1\ar j,3\ar l,1\ar q,2\ar s) $.
The last equation is due to the fact that $A_3$ is anti-symmetric.
$A_2(q,s).X(2\ar j,3\ar l,1\ar q,r\ar s)=A_2(q,s).X(1\ar j,3\ar l,2\ar q,r\ar s)$ for all $r>2$.
$A_2(q,s).X(2\ar j,3\ar l,2\ar q,r\ar s)=A_2(q,s).X(1\ar j,3\ar l,1\ar q,r\ar s)$ for all $r>2$.

Plugging the right hand side expressions into the expression of $C(2,3)_{jl}$ we get
$C(2,3)_{jl} = -C(1,3)_{jl}$. Hence $C(2,3) = -A_2.2.n.(n-2)!$. Similarly we can show that 
$C(2,k)=-A_2.2n(n-2)!$ for all $k\geq 3$.


Case $i>2,k>2$: Finally we have to show that $C(ik) = 0$ for $2<i$ and $2<k$.

If $i$ and $k$ are both greater than $2$, we have
 $A_1(q,s).X(i\ar j,k\ar l,1\ar q,1\ar s) = A_1(q,s).X(i\ar j,k\ar 2,1\ar q,2\ar s)$
and $A_2(q,s).X(i\ar j,k\ar l,1\ar q,r\ar s) = A_2(q,s).X(i\ar j,k\ar l,2\ar q,r\ar s)$.
Further, $A_3$ is anti-symmetric so 
$\sum_{q,s}A_3(q,s).X(i\ar j,k\ar l,1\ar q,2\ar s) = \sum_{q<s}(A_3(q,s)+ A_3(s,q)).
X(i\ar j,k\ar l,1\ar q,2\ar s) = 0$. Thus $C(i,k)=0 = 0.2n(n-2)!$ for all $i > 2, k>2$.

All these results combine to show that $v(E)$ is an eigenvector of $B$ with eigenvalue $2n(n-2)!$.

To see that there are several linearly independent vectors which are eigen vectors for the same
eigenvalue consider the following modifications in $E$. Begin with the observation that in blocks
$A_1,A_2$, and $A_3$ indices $1$ and $n$ have special role, while all other indices are equivalent.
That is, for each $k\in \{1,2,3\}$, $A_k(q,s) = A_k(q',s')$ where $q'=q$ if $q\in \{1,n\}$ otherwise 
$q'$ is any member of $\{2,\dots,n-1\}$ and $s'=s$ if $s\in \{1,n\}$ otherwise $s'$ is any member 
of $\{2,\dots,n-1\}$. So defining $A'_1,A'_2$ and $A'_3$ by exchanging the role of $n$ by any other 
index $\alpha$ in $\{2,\dots,n\}$ we again get an eigenvector with the same eigenvalue. These $n-1$ 
eigenvectors are linearly independent because in each vector the entry $B_3(2,n')$ is non-zero 
in exactly one vector, for each $n'\in \{2,\dots,n\}$.

There are additional $n-2$ similar sets of eigenvectors. Exchange the role of $2$ in $E$ by any 
index $\beta$ in the range $2,3,\dots,n-1$, see matrix $E'$. Once again we get a set of $n-1$ 
eigenvectors. Observe that among the diagonal blocks $22$ to $nn$ exactly one block is non-zero 
in each set. Hence these $(n-1)^2$ vectors are linearly independent.

\[E'= \left[ 
\begin{array}{cccccccc}
A_1   &  A_2   & \dots   & A_2   & A_3   & A_2   & \dots & A_2\\
A_2^T &  0     & \dots   & 0     &-A_2^T & 0     & \dots & 0  \\
\dots &  \dots & \dots   & \dots &\dots  & \dots & \dots & \dots\\
A_2^T &  0     & \dots   & 0     &-A_2^T & 0     & \dots & 0  \\ 
A_3^T &   -A_2 & \dots   & -A_2  & -A_1  & -A_2  & \dots & -A_2\\
A_2^T &  0     & \dots   & 0     &-A_2^T & 0     & \dots & 0  \\
\dots &  \dots & \dots   & \dots &\dots  & \dots & \dots & \dots\\
A_2^T &  0     & \dots   & 0     &-A_2^T & 0     & \dots & 0  \end{array}
\right ]\]

\begin{theorem} The dimension of the eigenspace of $B$ corresponding to eigenvalue $2n(n-2)!$ 
is at least $(n-1)^2$.
\end{theorem}

\end{document}